\documentclass{acmart}

\newcommand{\thm}[2]{\begin{#1} #2 \end{#1}}
\newcommand{\excess}{\mathrm{excess\:}}
\newcommand{\diameter}{\mathrm{diameter\:}}
\newcommand{\Sim}{\mathrm{Sim\:}}

\newtheorem{theorem}{Theorem}[section]

\begin{document}


\title{Combinatorial optimization in geometry\footnote{This is a
slightly modified version of a preprint first distributed in September
1996}}

%
\author{Igor Rivin Mathematics Institute, Warwick University\\
Coventry, The United Kingdom\\
and
\\
Mathematics Department, California Institute of   Technology
\\Pasadena, CA 91125\\
current address:\\
 Institut des Hautes \'Etudes Scientifiques, \\
35 Route de Chartres, 91440, Bures-sur-Yvette, France
}

\begin{abstract}
In this paper we extend and unify the results of \cite{steiner} and
\cite{rvol}. As a consequence, the results of \cite{steiner} are
generalized from the framework of ideal polyhedra in ${\bf H}^3$ to
that of singular Euclidean structures on surfaces, possibly with an
infinite number of singularities (by contrast, the results of
\cite{steiner} can be viewed as applying to the case of non-singular
structures on the disk, with a finite number of distinguished
points). This leads to a fairly complete 
understanding of the moduli space of such Euclidean structures and
thus, by the results of Epstein, Penner, (\cite{penner,ep,np}) the
author \cite{rvol,rivid}, and others, further
insights into the geometry and topology of the Riemann moduli space. 

The basic objects studied are the canonical \emph{Delaunay}
triangulations associated to the aforementioned Euclidean structures.

The basic tools, in addition to the results of \cite{rvol} and
combinatorial geometry are methods of combinatorial optimization --
linear programming and network flow analysis; hence the results
mentioned above are not only \emph{effective} but also
\emph{efficient}. Some applications of these
methods to three-dimensional topology are also given (to prove a
result of Casson's).
\end{abstract}

\maketitle

\section*{Introduction}

In the paper \cite{steiner} we gave the following description of the
angles of ideal polyhedra in $\mathbf{H}^3$: let $P$ be a combinatorial
polyhedron, and let $A: E(P) \rightarrow [0, \pi)$ be a
function. Then, there exists an ideal polyhedron combinatorially
equivalent to $P$, such that the exterior angle at every edge $e$ is
given by $A(e)$ if and only if the sum of $A(e)$ over all edges
adjacent to a vertex of $P$ is equal to $2\pi,$ while the sum of
$A(e)$ over any \textit{nontrivial cutset} of edges of $P$ (that is, a
collection of edges which separates the $1$-skeleton of $P$, but which
are not all adjacent to the same vertex) is strictly greater than
$2\pi.$ Furthermore, it was shown in \cite{rvol} that the dihedral
angles determine the ideal polyhedron up to congruence.

It was observed in \cite{rivid,rvol} that this was a
special case of the problem of characterizing the  \textit{Delaunay
tessellations} of singular Euclidean surfaces -- there is a canonical
way to associate ideal polyhedra to Delaunay triangulations of a
convex flat disk with convex polygonal boundary. The general situation
is described in detail below, but one of the goals of this paper is to
extend the characterization above to the completely general case of
singular Euclidean surfaces with boundary. This is duly done, see, eg,
Theorems \ref{mainthm}, \ref{versimple}, \ref{andreevthm}:

Consider a surface $S$, equipped with a Euclidean (or, more generally,
a similarity) structure $E$, possibly
with cone singularities. Assume that there is a discrete collection
$P=\{p_1, \ldots, p_n\}$ of distinguished points on $S$, and assume
that $P$ contains the cone points of $S$. There is a canonical
tessellation attached to the  triple $(S, E, P)$ -- the so-called
Delaunay tessellation (see, \textit{eg}, \cite{be,rvol}). The moduli
space $\mathcal{M}$ of such triples is then naturally decomposed into
disjoint subsets $\mathcal{M}_T$, corresponding to the different
combinatorial 
types $T$ of the Delaunay tessellation. This is a canonical
decomposition of $\mathcal{M}$. In the paper \cite{rvol} I studied the 
subsets $\mathcal{M}_T$, and showed that the \emph{dihedral angles} of the
Delaunay triangulation are natural coordinates (\emph{moduli}) for
$\mathcal{M}_T$, which induce on $\mathcal{M}_T$ the structure of a
convex polytope. 

The aformentioned decompositon of moduli space then becomes a
polyhedral complex, the top-dimensional cells of which corresponds to
Delaunay tessellations 
which are triangulations, while pairs of adjacent top-dimensional
cells differ combinatorially by a diagonal flip. This decomposition is
closely related to the well-known Harer complex (see, eg,
\cite{harer}). As mentioned above, the top-dimensional cells of this
complex are identified along some of their lower-dimensional faces,
while other lower dimensional faces correspond to degenerations of the
Euclidean structures of $(S, E, P)$. It is then clear that the
polyhedral structure of the cells $\mathcal{M}_T$ is of considerable
interest. However, in \cite{rvol} only an indirect description was
given -- $\mathcal{M}_T$ was shown to be a convex polytope by virtue
of being an image of another convex polytope under a fairly
complicated linear map. 
 The methods of \cite{steiner} come from hyperbolic
geometry and are based on the study of dihedral angles of
\emph{compact} hyperbolic polyhedra in \cite{rh}, so do not easily
generalize to the case of general singular Euclidean and similarity
structures alluded to above. In the current paper, 
methods of mathematical programming and the results of
\cite{rvol} are used to give a completely general extension of the
result of \cite{steiner} (described in the beginning of this
Introduction)  to Delaunay triangulations of arbitrary
singular surfaces (Theorems \ref{mainthm}, \ref{versimple},
\ref{andreevthm}). Since the arguments do not depend on the results of
\cite{steiner}, we have a different, essentially
combinatorial, proof of a principal result (Theorem 0.1) of that
paper (Theorem \ref{andreevthm} here). The other result of
\cite{steiner} -- the characterization of finite-volume polyhedra --
is, seemingly, not accessible by current methods. The above-mentioned
result permit us to get good understanding of the boundary structure
of the cells $\mathcal{M}_T$, and, consequently, of $\mathcal{M}$
itself. Since $\mathcal{M}_T$ fibers (in multiple ways) over the
moduli space of finite area hyperbolic structures on $S$, some
information is obtained about the latter moduli space.

In addition, the methods, combined with a geometric estimate, allow us
to give a description of dihedral angles of Delaunay tessellations of
$(S, E, P)$, where $(S, P)$ is not necessarily of finite topological
type (Theorem \ref{infmain}). This stops well short of solving the
moduli problem, unlike in the finite case, but a conjectural picture
seems fairly clear. 

The methods are also brought to bear onto some questions in
combinatorial geometry, and to provide efficient algorithms for
solving the ``inverse problem'' of determining when a combinatorial
complex, or a combinatorial complex equipped with dihedral angle data,
can be realized as the Delaunay tessellation of a singular Euclidean
surface. 

In addition, we use our methods to prove some observations of Casson
on ideally triangulated $3$-manifolds. That subject is not so far
removed from the geometry of similarity and Euclidean structures on
surfaces. Indeed, the basic idea of \cite{rvol} is to study the
similarity and Euclidean structures by means of constructing a
canonical hyperbolic polyhedral complex as a ``cone'' over the surface
being studied.

The plan of the paper is as follows: In section \ref{oldstuff}
the relevant definitions and results of \cite{rvol} are recalled. 
In section \ref{necessary} we describe a set of constraints which must
be satisfied by the dihedral angles of any (not necessarily Delaunay)
triangulation. In section \ref{suff} we show that these constraints
are actually sufficient under the assumption that the triangulation
is Delaunay, and refine them to a minimal set of constraints.  In
section \ref{ideal} we show how the results apply to 
ideal polyhedra, and in particular to characterize
infinite ideal polyhedra. In section \ref{moduli} we comment on the
boundary structure of the moduli space of singular Euclidean
structures, and describe a correspondence between the Euclidean and
hyperbolic structures, which hopefully clarifies the picture. 
In section \ref{lhyp} we
apply the methods of section \ref{suff} to the study of ideal
triangulations of $3$-manifolds. In section \ref{nflow} we give a
network flow interpretation of the results of section \ref{dihsec}. In
addition to the intrinsic interest, this allows us to give efficient
algorithms for deciding whether a weighed graph is the 1-skeleton of a
Delaunay triangulation (with weights being the dihedral angles). These
computational issues are discussed in section \ref{complexsec}. In
section \ref{someapps} we discuss some combinatorial-geometric
applications of the results of section \ref{dihsec}.

\section{Background}
\label{oldstuff}
\subsection{Singular similarity structures.} 
Consider an oriented surface $S$, possibly with boundary, and with a
number of distinguished points $\{p_1, \ldots, p_n\}$. A
\emph{similarity structure} on $S$ is given by an atlas for $S
\backslash \{p_1, \ldots, p_n\}$, such that the transition maps are
Euclidean similarities. A similarity structure induces a
holonomy representation $H_s$ of $\Gamma = \pi_1(S \backslash \{p_1, \ldots,
p_n\})$ into the similarity $\Sim (\mathbf{E}^2)\simeq
\widetilde{\mathbf{C}^*}$, where the tilde indicates the universal
cover. We define the  
\emph{dilatational   holonomy} as the induced representation $H_d:
\Gamma \rightarrow \mathbf{R}$, where
$H_d(\gamma) = \log \mbox{dilatation} H_s(\gamma) = \log |H_s(\gamma)|$. The
\emph{rotational holonomy} can almost be defined as $H_r(\gamma) = \arg
H_s(\gamma)$, but for the slight complication that we need to take the
argument in $\widetilde{\mathbf{C}^*}$, since we want to distinguish, eg, the
angle of $4\pi$ from one of $2\pi.$ This notion of argument is what is
used in the sequel. In particular, if $\gamma$ is a
loop surrounding one of the distinguished points $p_i$, then 
$H_r(\gamma)$ is the \emph{cone angle} at $p_i$. In the case where
$H_d(\Gamma) = \{0\}$, the similarity  structure is a \emph{singular
  Euclidean structure}, with cone angles defined by $H_r$ as above. In
the sequel, references to the \emph{holonomy} of a similarity
structure will actually mean the dilatational holonomy $H_d$. A
more concrete way to think of both similarity and Euclidean structure
is by assembling our surface out of Euclidean triangles, in the
pattern given by some complex $T$. In the case
where the lengths of the edges of the triangles being glued together
agree, then we have a singular Euclidean structure. If not, then we
have a similarity structure. In either case, the vertices of $T$ are
potentially cone points, with cone angles given by the sums of the
appropriate angles of the incident triangles.

Consider an oriented surface $S$, possibly with boundary, and assume 
that $S$ has a Euclidean metric (or, more generally, a similarity
structure -- with the exception of the results of Section \ref{ideal} ,
some of which depend on the metric Theorem \ref{infty}, the metric
structure or lack thereof plays no role in the arguments).
with cone singularities. We will be
dealing with \emph{semi-simplicial} triangulations of $S$, that is,
triangulations where two and where distinct closed cells might
intersect in a 
collection of lower-dimensional cells. All triangulations will be
assumed semi-simplicial, unless specified otherwise.  A subcomplex
$\mathcal{F}$ of $T$ will be called \emph{closed} if whenever an open
face $F$ is in $\mathcal{F}$, so are all of the faces of $\partial
F$. 

Assume now that the surface $S$ is equipped with a finite geodesic
semi-simplicial triangulation $T$, 
 such that the $0$-skeleton of $T$, which is denoted by $V(t)$,
contains all of the cone points of $S$. Each face of $T$ is then a
Euclidean triangle. There are two kinds of edges of $T$ -- the
interior edges, incident to two faces of $T$, and the boundary
edges, incident to only one face of $T$. 

\thm{definition}
{
\label{dihdef}
Let $e$ be an edge of $T$. First, suppose that
$e$ is a boundary edge, and let $t=ABC$ be a face of $T$ incident
to $e$, so that $e=AB$. Then the \emph{dihedral angle} $\delta(e)$ 
at $e$ is
the angle of $t$ at the vertex $C$. Now, assume that $e$ is an
internal edge of $T,$ so that $e$ is incident to $t_1 = ABC$ and 
$t_2 = ABD$, so that $e = AB.$ Then the dihedral angle at $e$ is
the sum of the angle of $t_1$ at $C$, and the angle of $t_2$ at
$D.$ The \emph{exterior} dihedral angle $\delta^\prime(e)$ at $e$ is
defined to be $\delta^\prime(e) = \pi - \delta(e)$. 

The \emph{cone angle} at an interior vertex $v$ of $T$ is the sum of
all the angles of the faces of $T$ incident to $v$ at $v$; the
\emph{boundary angle} at a boundary angle is defined in the same way. 
}

The following is a slight extension of \cite[Lemma 4.2]{rvol}.

\thm{lemma}
{
\label{cona}
The cone angle at an interior vertex $v\in V(T)$ is equal to the sum of the
exterior dihedral angles at the edges of $T$ incident to $v.$ At a
boundary vertex, the boundary angle is equal to the sum of the
exterior dihedral angles as above, less $\pi$.
}

\thm{proof}
{
First, let $v$ be an interior vertex.
Suppose that there are $n$ triangles $t_1, \dots, t_n$ 
incident to
$v$. The sum of all of their angles is $n \pi.$ The cone angle
at $v$ is the sum of the angles of the triangles $t_i$ at $v.$
On the other hand, the sum of the dihedral angles of the edges
$e_1, \dots, e_n$ 
incident to $v$ is the sum of all of the angles of $t_i$ \emph{not}
incident to $v$. Thus, 
\begin{equation}
\label{anglesum1}
n \pi = \mbox{Cone angle at $v$} + \sum_{i=1}^n \delta(e_i).
\end{equation}
The result follows by rearranging the terms.

If $v$ is a boundary vertex, and there are $n$ faces incident to $v$,
then there are $n+1$ edges, and equation \ref{anglesum1} becomes:
\begin{equation}
\label{anglesum2}
n \pi = \mbox{boundary angle at $v$} + \sum_{i=1}^{n+1} \delta(e_i),
\end{equation}
and the result follows by rearranging terms, as above.
}

\thm{observation} 
{
\label{gaussbonnet}
The sum of all of the dihedral angles
of all of the edges of $T$ is equal to the $\pi |V(T)|$ -- in 
combination with the Lemma above this gives the Gauss-Bonnet theorem
in this polyhedral context, since the curvature at an interior vertex 
$v$ of $T$ is defined to be $2\pi - \mbox{Cone angle at $v$}$, while 
the curvature at a boundary vertex is defined to be $\pi -
\mbox{boundary angle at $v$}$. 
}

\thm{proof}
{
Every angle of every face of $T$ is opposite to exactly
one edge of $T$.
}

The next theorem is \cite[Theorem 6.16]{rvol}.
\thm{theorem}
{
\label{uniqueness}
Let $\Delta: E(T) \rightarrow (0, 2\pi)$ be an assignment of dihedral
angles to the edges of $T$, and $H_d$ a holonomy representation.
There exists at most one singular similarity structure 
on $S$ with holonomy $H_d$ (and in particular,  at most one
singular Euclidean structure, up to scaling),
such that  so that $\delta(e) = \Delta(e),$ for every edge
$e\in E(T)$.
}

\thm{definition}
{
A triangulation $T$ with $\delta(e) \leq \pi$
for every interior $e\in E(T)$ is called a \emph{Delaunay triangulation}.
}

Let $\Delta$ be a map --  $\Delta: E(T) \rightarrow (0, 2\pi).$ When
does there exist a singular Euclidean metric on $S$ with the
dihedral angles prescribed by $\Delta$? It is clear that there are
certain linear constraints which must be satisfied -- to wit,
for every face $t=ABC$ of $T$, we must be able to find angles
$\alpha$, $\beta$, and $\gamma$, such that:

\begin{description}

\item[Positivity] All angles are strictly positive.

\item[Euclidean Faces] For every face $t$, $\alpha_t + \beta_t +
  \gamma_t = \pi.$ 

\item[Boundary edge dihedral angles] For every boundary edge $e = AB,$
  incident to the triangle $ABC$, $\gamma = \Delta(e).$

\item[Interior edge conditions] For every interior edge $e = AB,$
  incident to  triangles $ABC$ and $ABD$, $\alpha+\delta = \Delta(e).$
\end{description}

The system of equations and inequalities above specify a \emph{linear
program} $L$. The feasible region  
(set of solutions) of $L$ must be non-empty in order for us to have
any hope of having a singular Euclidean metric on $S$ with
prescribed dihedral angles. One of the principal results (Theorem 6.1) 
of \cite{rvol} is that if all of the dihedral angles are no greater
than $\pi$ (that is, the triangulation is Delaunay), then Conditions
1-4 are also sufficient, thus:

\thm{theorem}
{
\label{exist}
If the feasible region of $L$ is non-empty and every dihedral angle
is at most $\pi$, then there exists a similarity structure with any
prescribed holonomy $H_d$  and, in particular, a singular Euclidean metric on
$S$ with the prescribed dihedral angles -- this structure is unique
by Theorem \ref{uniqueness} (and the metric is unique up to scaling).
}

\section{Necessary conditions on dihedral angles}
\label{necessary}

In order for the linear program $L$ to have any chance of having
a solution, the dihedral angles $\Delta$ must satisfy some constraints.
Indeed, suppose $L$ has a solution. 
\thm{condition}
{
\label{dihpos}
All of the dihedral angles must be positive. 
}
Furthermore, the sum of all of the
dihedral angles of all of the edges of $T$ must be equal to the
sum of all of the angles of all of the faces of $T,$ or $F(T) \pi$.
On the other
hand, it is almost equally obvious that if $t_1, \dots, t_n \in F(T)$
is is some proper subset of the faces of $T$, then the sum of the
dihedral angles at the edges incident to one of the $t_i$ must be
strictly greater than the $n \pi$. In other words:

\thm{condition}
{
\label{sumdih}
Let $\mathcal{F} \subseteq F(T)$, and let $E(\mathcal{F})$ be the set
of all edges incident to an element of $\mathcal{F}$. Then  
\begin{equation}
\label{sumdihineq}
\sum_{e\in E(\mathcal{F})} \Delta(e)\geq \pi |\mathcal{F}|,
\end{equation}
with equality if and only if $\mathcal{F} = F(T)$ or
$E(\mathcal{F})=\emptyset.$
}

\thm{definition}
{
For a subcomplex $\mathcal{F}$ of $T$, we define the \emph{excess} of
$\mathcal{F}$ to be 
$$
\excess \mathcal{F} = \sum_{e\in E(\mathcal{F})} \Delta(e) - \pi |\mathcal{F}|,
$$
}

\thm{definition}
{
An edge $e$ of a subcomplex $\mathcal{F}$ of $T$ is a \emph{relative
boundary edge} of $\mathcal{F}$ if it is incident to a
top-dimensional face of $\mathcal{F}$ on exactly one side, and is not
a boundary edge of $T$.
}

The following lemma will prove useful in the sequel.

\thm{lemma}
{
\label{boundaryineq}
Let $\Delta: E(T) \rightarrow (0, \pi ]$ be an assignment of dihedral
angles to edges of $T$ satisfying condition \ref{sumdih}. 
Let $t_1, t_2, \ldots, t_n$ be a collection of closed triangles of $T$, and
let $\mathcal{F} = t_1 \cup t_2 \cup \ldots \cup t_n$ -- assume
$\mathcal{F} \neq T$. Let
$\sum_\partial \mathcal{F}$ be the sum of the dihedral angles of the
boundary edges of $\mathcal{F}$. Then
\begin{equation}
\label{doublebd}
0 < \excess \mathcal{F} < \sum_\partial \mathcal{F}.
\end{equation}
}

\thm{proof}
{
The first inequality of \ref{doublebd} is a restatement of Condition
\ref{sumdih}, applied to the subcomplex $\mathcal{F}$. To show the
second inequality, let 
$$
\overline{\mathcal{F}} = \bigcup_{t \in F(T) \backslash \{t_1,
\ldots, t_n\}} t.
$$
Evidently, 
$$
\mathcal{F} \cup \overline{\mathcal{F}} = T,
$$
while
$$
\mathcal{F} \cap \overline{\mathcal{F}} = \partial\mathcal{F}.
$$
Applying the conditions \ref{sumdih} to $\mathcal{F}$, we see that
\begin{equation}
\label{ineqf}
0 < \sum_{e\in E(\mathcal{F})} \Delta(e) -  \pi |F(\mathcal{F})|,
\end{equation}
while applying them to $\overline{\mathcal{F}}$, we see that 
\begin{equation}
\label{ineqf2}
0 < \sum_{e\in E(\overline{\mathcal{F}})} \Delta(e) -  \pi
|F(\overline{\mathcal{F}})|. 
\end{equation}

Note now that 
\begin{equation}
\label{ufaces}
|F(\mathcal{F})| + |F(\overline{\mathcal{F}})| = |F(T)|,
\end{equation}
while
\begin{equation}
\label{uedges}
\sum_{e\in E(\mathcal{F})} \Delta(e)+ \sum_{e\in
  E(\overline{\mathcal{F}})} \Delta(e) - \sum_\partial \mathcal{F} =
\sum_{e\in E(T)} \Delta(e).
\end{equation}
Adding inequalities \ref{ineqf} and \ref{ineqf2}, and applying equations
\ref{ufaces} and \ref{uedges}, we obtain
\begin{eqnarray*}
 & &  \sum_{e\in E(\mathcal{F})} \Delta(e) -  \pi |F(\mathcal{F})| \\
 & < & \left(\sum_{e\in E(\mathcal{F})} \Delta(e) -  \pi |F(\mathcal{F})|\right) +
         \left(\sum_{e\in E(\overline{\mathcal{F}})} \Delta(e) -
         \pi|F(\overline{\mathcal{F}})|\right) \\
  & = & \left(\sum_{e\in E(\mathcal{F})} \Delta(e) + \sum_{e\in
    E(\overline{\mathcal{F}})} \Delta(e)\right) - \left(\pi |F(\mathcal{F})| +
         \pi|F(\overline{\mathcal{F}})|\right) \\
  & = & \sum_{e\in E(T)} \Delta(e) +  \sum_\partial \mathcal{F} + \pi
  |F(T)|  =  \sum_\partial\mathcal{F},
\end{eqnarray*}
where the last equality is obtained by applying the equality case of
Condition \ref{sumdih}. 
}

\section{Sufficiency of conditions 2.1 and 2.2}
\label{suff}

Somewhat surprisingly, the trivial conditions \ref{dihpos} and \ref{sumdih}
guarantees that the linear program $L$ has a solution:

\thm{theorem}
{
\label{mainthm}
Let $T$ be a semi-simplicial triangulation of a surface $S$, possibly
with boundary, let $\Delta:E(T)\rightarrow (0, \pi]$ be a function
on the edges of $T$, $H_d$ a representation $\pi_1 (S \backslash V(T))
\rightarrow \mathbf{R}$. There exists a similarity structure with
holonomy $H_d$ on $S$, with the Delaunay triangulation of $(S,
\mathcal{E})$ combinatorially equivalent to $T$, and with dihedral
angles given by $\Delta$ if and only if conditions \ref{dihpos} and
\ref{sumdih} are satisfied.
}

We will use the Duality Theorem of linear programming to prove Theorem
\ref{mainthm}. Before stating the Duality Theorem we need to recall
some notions: A \emph{linear program} $L$ consists of a collection
$\mathcal{C}(L)$ of \emph{constraints}, which are linear equations or
(nonstrict) inequalities, and the \emph{objective function} $F(L)$,
which is a linear function. The set of points in $\mathbf{R}^n$
satisfying the constraints $\mathcal{C}(L)$ is called the
\emph{feasible region} of $L$, which is, by definition, a polyhedral
region, possibly unbounded, possibly not of full dimension, and
possibly empty. The \emph{solution} of the linear program is a point
where the objective function attains an extremum (which may be a maximum or a
minimum). The value of the objective function at the solution is the
\emph{objective} of $L$. If the feasible region of $L$ is empty, the
program is said to be infeasible. Now we can state the Duality Theorem:

\begin{theorem}[Duality Theorem of Linear Programming]
\label{duality}
Let $P$ be a linear program of the form:

Minimize ${\bf c}^\perp {\bf x}$, subject to the constraints
${\bf A} {\bf x} = {\bf a}$, ${\bf x} \geq 0$.

Then the dual  of $P$ is the program $P^*$:

Maximize ${\bf a}^\perp {\bf \lambda}$, subject to the constraints
${\bf \lambda}^\perp {\bf A} \leq {\bf c}$.

The feasible region of $P$ is nonempty if and only if the 
objective of $P^*$ is bounded. Conversely, the feasible region
of $P^*$ is nonempty if and only if the objective of $P$ is bounded.
If neither feasible region is empty, then the values of the objective
functions of $P$ and $P^*$ are equal.
\end{theorem}

\medskip\noindent
\textit{Remark.}
While the \emph{primal} program $P$ in the statement of Theorem
\ref{duality} appears to be of a very special form, it is not hard to
see that any linear program can be written in this form. Indeed, if
our linear program asked us to maximize the objective, we can always
convert it to a minimization problem by multiplying $\mathbf{c}$ by
$-1$. If our program did not require the variables to be non-negative,
we can always replace a variable $x$ by $x_+ - x_-$, where $x_+$ and
$x_-$ are both required to be non-negative. If the program had some
inequalities of the form $\mathbf{a_i}^\perp \mathbf{x} \geq a$, or of
the form $\mathbf{a_i}^\perp \mathbf{x} \leq a$ we can first convert
all the inequalities of the first type to those of the second type by
negation, and then convert them to equations by introducing
\emph{slack variables} $x_i \geq 0$, and requiring 
$\mathbf{a_i}^\perp \mathbf{x} + x_i = a.$ Similarly, any program can
be made to look like the dual program $P^*$ in the statement of
Theorem \ref{duality}.

Program $L$ is almost in the primal form needed by the duality
theorem, but for two differences: there is no 
objective function, and we want the primal variables (the angles
of the triangles) to be strictly positive, rather than just non-negative.

For the moment, let us sweep these issues under the rug, by 
setting the objective function to be $0$, and allowing the
angles to vanish -- it will be quickly apparent how to fix
things up later. Let the modified program be $L_1$. The dual 
$L_1^*$ of $L_1$ is the following:

\begin{description}

\item[The dual program]
Maximize $F({\bf u}, {\bf v}):\pi \sum_{t \in F(T)} u_t + \sum_{e\in E(T)} \Delta(e)v_e$,
subject to the conditions
$u_t+v_e\leq 0$ whenever $e$ is an edge of $t$.

\end{description}

\thm{theorem}
{
\label{maindual}
Assume that the conditions described in conditions \ref{dihpos} and
\ref{sumdih} are satisfied.
Then the objective function of $L_1^*$ is nonpositive. It equals 
$0$ if and only if there is a $u$, such that $u_t=-v_e=u$, for
all $t\in F(T)$, $e\in E(T)$.
}

\thm{proof}
{
\begin{description}
\item[Observation 1] Note that if there is a $u$ as required in the
statement of  the Theorem, the objective function is, indeed, equal to
$0$. This is nothing other than the equality case of Condition \ref{sumdih}.
\end{description}

Now, let $u = \min(u_1, \dots, u_{F(T)})$. 
Let $u_i^{(1)} = u_i - u$, and let $v_j^{(1)} = v_j + u$, for
all values of the indices. The new variables are still feasible
for $L_1^*$, and by the observation 1 above, this transformation
does not change the value of the objective. Furthermore, if
$u_i^{(1)} = 0$ for all $i$, then the objective is non-positive,
and is equal to zero only if all of $v_j^{(1)}$ are equal to zero 
as well -- this is so, since all of the $v_j$ must be non-positive, and
all of their coefficients are positive, by Observation \ref{dihpos}. 
Assume, then, that $u_t^{(1)}>0$ for $t \in \mathcal{F}^{(1)}$; 
$\mathcal{F}^{(1)}$ is a \emph{proper} subset of $F(T)$ by
construction.

\begin{description}
\item[Observation 2] Suppose now that $u_t^{(1)} = u$, for all 
$t \in \mathcal{F}_1$. Then 
$$F\leq \pi \sum_{t\in \mathcal{F}^{(1)}}u-
\sum_{e\in E(\mathcal{F}^{(1)})} \Delta(e) u < 0,$$ by
\ref{sumdihineq}.
\end{description}

Now, let $u^{(1)}  = \min_{t \in \mathcal{F}^{(1)}}(u_t^{(1)})$, and
let $u_t^{2}=u_t^{(1)}=0$ if $t \notin \mathcal{F}^{(1)}$, and otherwise
$u_t^{(2)} = u_t^{(1)}-u^{(1)}$. Likewise, $v_e^{(2)} = v_e^{(1)}$
if $e \notin E(\mathcal{F}^{(1)}$, and otherwise 
$v_e^{(2)}= v_e^{(1)} + u^{(1)}$. 

This still leaves us in the feasible region of $L_1^*$ and strictly
increases the value of the objective (by Observation 2). The new 
nonzero set $\mathcal{F}^{(2)}$ is a proper subset of $\mathcal{F}^{(1)}$,
and we can repeat this process. In the end, we will wind up with 
a feasible point ${\bf u}^{(k)}, {\bf v}^{(k)}$, with 
${\bf u}^{(k)} = {\bf 0}$,  where the value
of the objective is non-positive (by Observation 1), and strictly
greater than the value of the objective at ${\bf u}, {\bf v}$ (by
Observation 2), thus completing the proof.
}

The above theorem shows that the feasible region of $L_1$ is 
non-empty. In order to find a solution with strictly positive 
angles, we write our angles $\alpha_i$ as $\alpha_i = \beta_i+\epsilon$.
We require all of the $\beta_i$ to be non-negative, and
our new objective is simply $-\epsilon.$
Call the resulting program $L_2$. Its dual $L_2^*$ has the following
form:

\begin{description}
\item[Second dual]
Maximize $F({\bf u}, {\bf v}):\pi \sum_{t \in F(T)} u_t + \sum_{e\in E(T)} \Delta(e)v_e$,
subject to the conditions
$u_t+v_e\leq 0$ whenever $e$ is an edge of $t$, and also
$3  \sum_{t \in F(t)} u_t + 
2 \sum_{e\in E(T), e\notin \partial T} v_e + 
\sum_{e\in E(T), e \in \partial T} v_e \leq -1$.
\end{description}

\thm{theorem}
{
\label{modified}
The optimal value of the objective of $L_2^*$ is strictly negative.
}

\thm{proof}
{
Suppose the contrary. Then $u_t = -v_e = u$, for some $u$,
by Theorem \ref{maindual}. However, in that case the last inequality
of $L_2^*$ is not satisfied, since the left hand side of the last
constraint of $L_2^*$ vanishes. Indeed, it is equal to
$$
u \left(3  \sum_{t \in F(t)} 1 - 2 \sum_{e\in E(T), e\notin \partial
    T} 1 - \sum_{e\in E(T), e \in \partial T} 1 \right).
$$
Observe, however, that $3 \sum_{t \in F(t)} 1$ counts each
non-boundary edge with multiplicity $2$, and each boundary edge with
multiplicity $1.$
}

\section{Dihedral angles of Delaunay triangulations}
\label{dihsec}

In the preceding section it was shown that in order for the linear program
$L$ to have a solution, it is necessary and sufficient that the putative
dihedral angles $\delta(e)$ satisfy the inequalities \ref{dihpos} and
\ref{sumdih}. Below, we use these to derive another set of necessary conditions
\ref{kappaineq}, and show that under the additional assumption that
the original dihedral angles do not exceed $\pi$, these are equivalent
to the inequalities \ref{sumdih}. The virtue of the inequalities \ref{kappaineq}
is that it is easier to interpret them geometrically. A special case of them
is one of the main results of \cite{steiner} -- the connection will be explained
in section  \ref{ideal}.

First, some definitions (these are the analogs of the definitions of section
\ref{oldstuff} without reference to face angles) -- as before $T$ is a
triangulation of $S$, $V(T)$ is the set of vertices of $T$, $E(T)$ is
the set of edges, and $F(T)$ is the set of faces. We assume that each
edge $e\in E(T)$ is given a weight $\delta(e)\in \mathbf{R}$.  
The weight $\delta(e)$ will be called the dihedral angle at $e$. 
The quantity $\delta^\prime(e)=\pi - \delta(e)$ will be called the
\emph{exterior dihedral angle} at $e$. 

If $v$ is a vertex of $T$, then the \emph{cone angle} $C_v$
is defined to be $$C_v = \sum_{\mbox{$e$ incident to $v$}}
\delta^\prime(e),$$ while the \emph{curvature} $\kappa_v$ is defined
to be $\kappa_v = 2\pi - C_v.$ 

\medskip\noindent
\textit{Note.} The cone angle of a boundary point is, thus, not the
same as the boundary angle (in the language of section \ref{oldstuff},
but smaller by $\pi$.

Below, it will often be useful to talk of the \emph{Poincar\'e dual} of $T$. Recall that this
is the complex $T^*$, such that the set of vertices of $T^*$ is in one-to-one
correspondence with the set of faces of $T$, the set of edges of $T^*$ are in
one to one correspondence with the edges of $T$ -- two vertices of
$T^*$ are joined by an edges if and only if the corresponding faces of $T$ share
an edge, and finally the faces of $T^*$ correspond to the vertices of $T$ --
the vertices of a face $v^*$ of $T^*$ correspond to the faces of $T$ incident
to the corresponding vertex $v$. 

\thm{definition}
{
A subcomplex $\mathcal{F}$ of $T$ is \emph{closed}, 
if whenever a cell $t$ is in $\mathcal{F}$, then so are all of the
lower-dimensional cells incident to $t$.
}

\thm{definition}
{
The \emph{total curvature} of a subcomplex $\mathcal{F}$ of $T$ is
defined as
$$K(\mathcal{F}) = \sum_{v\in V(\mathcal{F})} \kappa_v.$$
}

\textit{Notation.}
Let $\mathcal{F}$ be a subcomplex of $T$, and let
$E^\prime(\mathcal{F})$ be the set of those edges $e$ of $T$ which are
not edges of faces of $\mathcal{T}$, but such that at least one
endpoint of $e$ belongs to $\mathcal{T}$. For each edge $e$ of $T$,
define $n_{\mathcal{F}}(e)$ to be the number of endpoints of $e$ which
belong to $\mathcal{F}$.

For example, if $e\in E(\mathcal{F})$, then $n_{\mathcal{F}}(e) = 2$;
if $e \notin E(\mathcal{F}) \cup E^\prime(\mathcal{F})$, then
$n_{\mathcal{F}}(e) = 0$.

\thm{theorem}
{
\label{kappaineq}
For every non-empty subcomplex $\mathcal{F}$ of $T$
the following are equivalent:

\begin{description}

\item[(a)]
$$
\sum_{e\in E^\prime(\mathcal{F})} n_{\mathcal{F}}(e)
\delta^\prime(e)\geq  2\pi \chi (\mathcal{F})- K(\mathcal{F}), 
$$
with equality if and only if $\mathcal{F} = T.$

\item[(b)] Conditions \ref{dihpos} and \ref{sumdih} hold.
\end{description}
}

\thm{proof}
{
We will show that (b) $\Rightarrow$ (a); the converse is immediate.
\begin{equation}
\label{kappasum}
K(\mathcal{F}) = \sum_{v\in V(\mathcal{F})} 
(2\pi - \sum_{\mbox{$v\in e$}} \delta^\prime(e))=
2\pi |V(\mathcal{F})| - \sum_{v\in V(\mathcal{F})}\sum_{\mbox{$v \in e$}}  
\delta^\prime(e).
\end{equation}
The last sum of equation \ref{kappasum} can  be rewritten thus:
\begin{equation}
\label{kappa2}
\sum_{v\in V(\mathcal{F})}\sum_{\mbox{$e$ incident to $v$}}
\delta^\prime(e) = 2 \sum_{e\in E(\mathcal{F})} \delta^\prime(e) +
\sum_{e\in E^\prime(\mathcal{F})}  n_{\mathcal{F}}(e) \delta^\prime
(e). 
\end{equation}

Finally, using the definition of $\delta^\prime(e)$, and combining equations
\ref{kappasum} and \ref{kappa2} it follows that
\begin{equation}
\label{kappaeq}
K(\mathcal{F}) = 2\pi \chi(\mathcal{F}) + 
2(\sum_{e\in E(\mathcal{F})} \delta(e) -\pi |F(\mathcal{F})|)
 - \sum_{e\in E^\prime(\mathcal{F})} \delta^\prime(e).
\end{equation}
Now, assume that the dihedral angles satisfy the inequalities \ref{sumdih}. This
means that the middle term on the right hand side of equation \ref{kappaeq} is
non-negative (and strictly positive unless $\mathcal{F} = T$).
}

\subsection{Some corollaries and refinements}
\label{corref}

Here are some easy consequences of Theorem \ref{kappaineq}:

\begin{description}

\item[Special case 1] $\mathcal{F} = T$. Then, since $E^\prime(\mathcal{F}) =
\emptyset$, 
we just get the Gauss-Bonnet theorem (with curvature defined in terms of the
dihedral angles).

\item[Special case 2]  The complement of $\mathcal{F}$ is an annulus
containing no vertices of $T$ -- this corresponds to a simple cycle in
the Poincar\'e dual $T^*$. For every edge in $E^\prime(\mathcal{F})$,
$n_e = 2$.   Theorem \ref{kappaineq} just says that $\sum_{e \in
E^\prime} \delta^\prime(e) > 0.$ 

\item[Special case 3] For every edge in $E^\prime(\mathcal{F})$, $n_e = 1$. This will
hold in the case where $\mathcal{F}$ is a subset of $T$ with a collar, and removing
the edges in $E^\prime$ separates the $1$-skeleton of $T$ into (at least) two 
connected components. In this case, Theorem \ref{kappaineq} tells us that
$$
\sum_{e\in E^\prime(\mathcal{F})} \delta^\prime(e)\geq 
2\pi \chi (\mathcal{F})- K(\mathcal{F}).
$$
\end{description}

A drawback of both conditions \ref{sumdih} and Theorem \ref{kappaineq}
is that they require their respective inequalities to be checked for
\emph{every} subcomplex of $T$ in order to verify whether a given
assignment of dihedral angles is admissible. It is not hard to see
that this requirement can be weakened somewhat.

Let $\mathcal{F}$ be the subcomplex in question.

\thm{observation}
{
\label{obs1}
It can assumed that every edge in $E(\mathcal{F})$ is an edge
of at least one face.
}

\thm{observation}
{
\label{obs2}
The $1$-skeleton of the Poincar\'e dual $\mathcal{F}^*$ of
$\mathcal{F}$ can be assumed connected -- this is somewhat 
stronger than saying that $\mathcal{F}$ is connected.
}

\thm{observation}
{
\label{obs3}
It can be assumed that the $1$-skeleton of $\overline{\mathcal{F}}^*$
-- the Poincar\'e dual of the complement of $\mathcal{F}$ -- has no
isolated vertices.  
}

It is straightforward to check all of the above observations.

In the special situation where $\delta(e)\leq \pi$ for all $e\in E(T)$ -- that is,
$\delta$ is Delaunay, one can further assume that every face $f$ of 
$\overline{\mathcal{F}}$ is adjacent to at most one face of $\mathcal{F}$. Otherwise,
adjoin $f$ to $\mathcal{F}$, to create a new complex $\mathcal{F}^\prime$. 
This complex has one more face than $\mathcal{F}$, but its sum of dihedral angles
is at most $\pi$ greater than that of $\mathcal{F}$. Hence, it is enough to check
that $\mathcal{F}^\prime$ satisfies the hypotheses of Theorem
\ref{kappaineq}, or the conditions \ref{sumdih}.

\thm{definition}
{
A \emph{simple} subcomplex of $T$ is a subcomplex $\mathcal{F}$ such
that both $\mathcal{F}$ and $T \backslash \mathcal{F}$ are connected.
}

\thm{theorem}
{
\label{versimple}
To verify that Conditions 2.1 and 2.2 (or equivalently, conditions of
Theorem \ref{kappaineq}[(b)]) hold for \emph{all} subcomplexes
of $T$, it is necessary and sufficient to check them for simple
subcomplexes. 
}

\thm{proof}
{
By Observation \ref{obs2}, it is enough to check the connected
subcomplexes of $T$. If such a subcomplex $\mathcal{F}$ is
simple, then we are done. Otherwise, its complement is not
connected. Let $\mathcal{C}$ be a connected component of the
complement of $\mathcal{F}$.
\thm{lemma}
{
The complex $\mathcal{C}$ is simple.
}

\thm{proof}
{
By construction, $\mathcal{C}$ is connected. Also, every point of 
$T \backslash \mathcal{C}$ can be connected by a path to a point of
$\mathcal{F}$. Since $\mathcal{F}$ is assumed connected, the lemma
follows. 
}

Consider $\mathcal{F}^\prime = \mathcal{F} \cup \mathcal{C}$. By
assumption, $\mathcal{F}^\prime \neq T.$ Now

\begin{eqnarray}
 & & \left(\sum_{e\in E(\mathcal{F})} \Delta(e) -  \pi
  |F(\mathcal{F})|\right) +
\left(\sum_{e\in E(\mathcal{C})} \Delta(e) -  \pi
  |F(\mathcal{C})|\right) - \sum_\partial \mathcal{C} \\
 & = &
\left(\sum_{e\in E(\mathcal{F}^\prime)} \Delta(e) -  \pi
  |F(\mathcal{F})^\prime|\right).
\end{eqnarray}

By lemma \ref{boundaryineq} (more precisely, a version for simple
complexes), 
\begin{equation}
\left(\sum_{e\in E(\mathcal{C})} \Delta(e) -  \pi
  |F(\mathcal{C})|\right) - \sum_\partial \mathcal{C} < 0, 
\end{equation}
thus,
\begin{equation}
\left(\sum_{e\in E(\mathcal{F})} \Delta(e) -  \pi
  |F(\mathcal{F})|\right) > 
\left(\sum_{e\in E(\mathcal{F}^\prime)} \Delta(e) -  \pi
  |F(\mathcal{F})^prime|\right).
\end{equation}

Thus, in order to check that condition \ref{sumdih} holds for
$\mathcal{F}$, it is enough to check that it holds for
$\mathcal{F}^\prime$. The proof of Theorem \ref{versimple} is finished
by the obvious induction 
argument on the number of connected components of the complement of
$\mathcal{F}$. 
}

In the case where  $T$ is a genuine simplicial complex (that is, two
cells intersect in a lower-dimensional cell), simple subcomplexes
corresponds to \emph{non-coterminous minimal cutsets} of edges of $T$:

\thm{definition}
{
A collection $\mathcal{C}$ of edges of $T$ is a \emph{cutset} if
removing the edges in $\mathcal{C}$ disconnects the $1$-skeleton of
$T$. A cutset $\mathcal{C}$ is $\emph{minimal}$ if no proper subset of
$\mathcal{C}$ 
is a cutset. A cutset $\mathcal{C}$ is \emph{coterminous} if all of
the edges in $\mathcal{C}$ are incident to the same vertex.
}

In other words, a minimal cutset corresponds to a separating simple
curve in the Poincar\'e dual $T^*$ of $T$ (the curve need not be closed if $T$
has boundary). A coterminous cutset corresponds to a boundary of a
face in the dual. The non-coterminous simple cutset corresponding to
the subcomplex $\mathcal{F}$ is nothing other than the set of edges
$E^\prime(\mathcal{F})$ define just before the statement of theorem
\ref{kappaineq}. 

In the case where the surface $S$ is a flat disk, theorems
\ref{maindual}, \ref{exist}, \ref{kappaineq}, and \ref{versimple}
immediately imply:

\thm{theorem}
{
\label{andreevthm}
Let $T$ be a triangulation of the disk, let $\Delta: E(T) \rightarrow
(0, \pi]$ be an assignment of dihedral angles to the edges of $T$, let
$V_\partial(T)$ be the set of boundary vertices of $T$ and
let $\Lambda: V_\partial(T) \rightarrow (0, \pi]$ be the assignment of
boundary angles. Then, there exists a collection of points $p_1,
\ldots, p_{V(T)}$ in the plane $\mathbf{E}^2$ whose Delaunay
triangulation is combinatorially equivalent to $T$, has dihedral
angles given by $\Delta$, and whose convex hull is a polygon with
angles given by $\Lambda$, if and only if:

\begin{itemize}

\item\label{sumvv} If $v$ is an interior vertex of $T$, then
\begin{equation}
\sum_{\mbox{$e \in E(T)$ incident to $v$}} (\pi - \Delta(e)) = 2\pi.
\end{equation}

\item\label{sumvb} If $v$ is a boundary vertex of $T$, then
\begin{equation}
\sum_{\mbox{$e \in E(T)$ incident to $v$}} (\pi - \Delta(e)) + (\pi -
\Lambda(v)) = 2\pi.
\end{equation}

\item\label{sumintcyc} If $\mathcal{E}$ is a non-coterminous minimal cutset
corresponding to a simple closed curve in the Poincar\'e dual $T^*$ of
$T$, then
\begin{equation}
\sum_{e \in \mathcal{E}} (\pi - \Delta(e)) > 2\pi.
\end{equation}

\item\label{sumextcyc} If $\mathcal{E}$ is a non-coterminous minimal
  cutset not corresponding to a simple closed curve in $T^*$, and
  $\mathcal{E}$ separates the boundary vertices of $T$ into two groups
  $v_1, \ldots, v_k$ and $v_{k+1}, \ldots, v_{|V_\partial(T)|}$, then
\begin{equation}
\sum_{e \in \mathcal{E}} (\pi - \Delta(e)) + \sum_{i=1}^k (\Pi-\Lambda(v_i))> 2\pi.
\end{equation}

\end{itemize}
}

\medskip\noindent
\textit{Remark.} Theorem \ref{andreevthm} is nothing but one of the
main results (Theorem 0.1) of \cite{steiner}, stated in a different
language (without mentioning convex ideal polyhedra), so we have
succeeded in deducing that theorem in a purely combinatorial way from
the results of \cite{rvol}. In addition, Theorem \ref{versimple} is
seen to be a direct generalization of \cite{steiner}[Theorem 0.1] to
a characterization of dihedral angles of Delaunay triangulations of
arbitrary, possibly singular, Euclidean surfaces.

\section{A network flow approach}
\label{nflow}

There is an alternative way to prove Theorem \ref{maindual} which uses the
Max Flow-Min Cut theorem of network flow instead of the duality 
theorem of linear programming. It should be noted that the difference 
between the two arguments is largely superficial, since the proof of 
Theorem \ref{maindual} can be seen to essentially prove the Max Flow-Min 
Cut theorem. There are two reasons to set up the question as a result 
on network flow. The first is that the proof (hopefully) becomes 
clearer and more intuitive. The second is that the special sorts of 
linear programs that arise in the theory of network flow have been 
heavily analyzed from the viewpoint of complexity, which will allow 
us to give a very satisfactory estimate (section \ref{complexsec}) 
for the running time of an algorithm to determine whether a structure 
with prescribed dihedral angles actually exists.

A \emph{network} is a directed (multi)graph $N$, with two distinguished 
vertices -- $s$ (the source), and $S$ (the sink). Each edge of $N$ 
has a certain \emph{capacity}, which is a real number, which is an 
upper bound on the amount of the commodity which can flow from the 
tail to the head of the edge. A \emph{cutset} $C$ of $N$ is a collection of 
edges, the removal of which leaves $s$ and $S$ in two different 
connected components of $N \backslash C$. The \emph{capacity of the
  cutset} 
$C$ is simply the sum of the capacities of the edges comprising $C$. 
The capacity of an empty collection of edges is, of course, $0$.

The Max Flow-Min Cut theorem of network flow (see, for example,
\cite[Chapter 7]{vlw}) says the following:

\begin{theorem}[Max Flow -- Min Cut]
\label{mfmc}
The maximal amount of a commodity that can flow from the source $s$ to 
the sink $S$ of a network $N$ is equal to the capacity of the smallest 
cutset in $N$.
\end{theorem}

This theorem can be proved in a number of ways. The interested reader 
can adapt the proof of Theorem \ref{maindual} to show Theorem
\ref{mfmc}. 

To use Theorem \ref{mfmc} for our purposes, we need to set up a 
network of a special sort, starting with a triangulation $\mathcal{T}$.

The vertices of this network are divided into four classes: $\{s\}$, 
$F(\mathcal{T})$, $E(\mathcal{T})$, and $\{S\}$.

The source $s$ is connected to all of the vertices corresponding to 
the faces of $\mathcal{T}$ with edges of capacity $1$. We call them
\emph{edges of level $1$}. 
Every vertex 
corresponding to a face $t$ is connected to the three vertices, corresponding 
to the three edges of $t$ by edges of capacity $1$. These are edges of 
level $2$. Finally, each $e\in 
E(\mathcal{T})$ is connected to the sink $S$ by  an edge of capacity 
$\delta(e)$. These are edges of level $3.$

The following statement is self-evident:

\thm{observation}
{
There exists a solution of the linear program $L_{1}$ if and only if 
the maximal flow through the above-constructed network $N_{1}$ is 
equal to $|F(\mathcal{T})|$.
}

\begin{proof}[of theorem \ref{mainthm}]
Consider a cut $C$ of $N_{1}$. This will have some edges at level $1$, 
removing which which will disconnect a subset $F_{0}$ of the faces 
of $\mathcal{T}$ from the source. It is then not necessary to remove any 
edges of level $2$ emanating from $F_{0}$. Let $F(\mathcal{T}) 
\backslash F_{0} = F_{1}$. Let $F_{2} \subseteq F_{1}$ be those 
faces $f$ for which the cutset contains all three edges of level $2$ 
emanating from $f$. Finally, $F_{3} = F_{1} \backslash F_{2}$. All 
of the edges of level $3$ (indirectly) emanating from $F_{3}$ must be 
in the cutset $C$. These are precisely the edges corresponding to the 
edges of the subcomplex of $\mathcal{T}$ whose faces are in $F_{3}$.

What is the capacity of $C$? Evidently, it is equal to

$$
\label{thesum}
F_{0}+ 3 F_{2}+ \sum_{e\in E(F_{3})} \delta(e).
$$
If we want the flow through $N_{1}$ to be $F$, we must have
$$
F_{0}+ 3 F_{2}+ \sum_{e\in E(F_{3})} \delta(e) \geq F,
$$
Or, noting that $F_{1} = F - F_{0}$, 
$$
3 F_{2}+ \sum_{e\in E(F_{3})} \delta(e) \geq F_{1}.
$$

In the special case where $F_{2} = 0$, it follows that
\begin{equation}
\label{geq}
\sum_{e\in E(F_{3})} \delta(e) \geq F.
\end{equation}

Since $F_{3}$ could have been any subset of $F$, it follows that 
condition \ref{sumdih} is necessary (this is, in any event, self evident).

On the other hand, if condition \ref{sumdih} holds for any subcomplex 
of $\mathcal{T}$, substituting the inequality \ref{sumdih} into 
\ref{thesum}, we see that 
$$c(C)\geq F_{0}+ 3 F_{2}+ F_{3} = F + 2 F_{2} \geq F.$$
\end{proof}

The above result is not quite what is required: the non-strict
inequality \ref{geq} implies the existence of a consistent assignment 
of angles to the faces of $\mathcal{T}$, but some of these angles may be 
equal to $0$. In order to have the angles strictly positive, we must 
modify the weights in the network $N_{1}$ as follows:

Modify the capacity of level $1$ edges to be $1-3\epsilon$; those of 
level $2$ to be $1-\epsilon$, and those of level $3$ to be 
$1-2\epsilon$. Call the resulting network $N_{2}$. The existence of 
a consistent assignment of angles to the faces of $\mathcal{T}$ where all 
the angles are no smaller than $\epsilon$ is obviously equivalent to 
the maximal flow in $N_{2}$ being equal to $F(\mathcal{T})(1-3 \epsilon)$.

Consider now a cut $C$ in $N_{2}$, with notation as before. The
capacity of $C$ will be: 

\begin{equation}
\label{epseq}
(1-3\epsilon)F_{0} + 3 (1-\epsilon)F_{2} + \sum_{e\in E(F_{3}} 
\delta(e) - 2\epsilon E(F_{3}).
\end{equation}

Assume that $F_{2} = 0$. The MinCut condition together with expression 
\ref{epseq} gives:
$$
(1-3\epsilon)F_{0} + \sum_{e\in E(F_{3}} 
\delta(e) - 2\epsilon E(F_{3}) \geq (1-3\epsilon) F,
$$
or
$$
\sum_{e\in E(F_{3}} \delta(e) \geq (1-3\epsilon)F_{3} + 2\epsilon 
E(F_{3}).
$$
If $\emptyset \subset F_{3}\subset F$, then $E(F_{3}) > 
\frac32 F_{3}$, and so
$\sum_{e\in E(F_{3}}\delta(e) > F_{3}.$

Suppose now that $\delta(F') - F' \geq \psi > 0$, for all proper 
non-empty subsets $F'$ of $F$.

Then the expression \ref{epseq} is no smaller than
\begin{equation}
\label{ineqq}
(1-3\epsilon)F_{0} + 3 (1-\epsilon)F_{2} + F_{3}+\psi- 2\epsilon E(F_{3}).
\end{equation}

The edges of $F_{3}$ can be divided into two classes: interior 
edges (those incident to two triangles of $F_{3}$) -- these number
$E_{i}(F_{3})$, and boundary edges of $F_{3}$ -- those incident to 
only one triangle. These number $E_{\partial}(F_{3})$. Since $F_{3}$ is 
a proper subset of $F$, $E_{\partial}(F_{3})>0$. Clearly,
\begin{equation}
2 E(F_{3}) = 3 F_{3} + E_{\partial}(F_{3}).
\end{equation}
The lower bound \ref{ineqq} can thus be rewritten as
\begin{eqnarray}
c(C) & \geq &(1-3\epsilon)F_{0} + 3 (1-\epsilon)F_{2} +
(1-3\epsilon)F_{3}+\psi- \epsilon E_{\partial}(F_{3}) \\
& = & (1-3\epsilon) F + \psi +  (1-3\epsilon)F_{2} - \epsilon
E_{\partial}(F_{3}). 
\end{eqnarray}
Since $\epsilon\leq\frac13$ (a triangle cannot have all angles greater 
than $\pi/3$), we get
\begin{equation}
c(C) -(1-3 \epsilon)F\geq \psi - \epsilon E_{\partial}(F_{3}).
\end{equation}

In other words, if $\sum_{e\in E(F')} \delta(e) - F'\geq \psi$ for 
every $\emptyset \subset F'\subset F$, there is a solution of 
the linear program $L_{1}$ with all face angles of all triangles no 
smaller than $\psi/E(\mathcal{T})$.
\qed

\section{Delaunay triangulations of infinite sets of points}
\label{ideal}
In this section we will show that theorem \ref{mainthm} for singular
Euclidean structures can be extended without change to infinite
locally finite complexes: 

\thm{theorem}
{
\label{infmain}
Let $T$ be an infinite but locally finite complex, and let $\Delta:
E(T) \rightarrow (0, \pi]$. Then there exists a singular euclidean
structure on $T$, with cone points at vertices of $T$, whose
Delaunay triangulation is combinatorially equivalent to $T$, and whose
dihedral angles are given by $\Delta$ if and only if each finite
subcomplex $\mathcal{F} \subset T$ with $E(\mathcal{F}) \neq
\emptyset$ has positive excess.
}

There are two ingredients in the argument. The first (Lemma
\ref{linref}) is an extension of section \ref{suff}, the second
(\ref{infty}) is a geometric estimate which will enable us to extract
the necessary subsequences.

\thm{lemma}
{
\label{linref}
Let $T$ be a complex, and let $\Delta: E(T) \rightarrow (0, \pi]$.
Then there exists a Euclidean structure with cone angles at vertices
of $T$, and dihedral angles given by $\Delta$, except at the boundary
edges of $T$, where they are smaller than prescribed by $\Delta$ if
the excess of any subcomplex $\mathcal{F}$ of $T$, such
that $E(\mathcal{F}) \neq \emptyset$ is positive.
}

\medskip\noindent
\textit{Remark.} Lemma \ref{linref} can be viewed as a relative
version of theorem \ref{mainthm}.

\thm{proof}
{
The argument parallels very closely that of section \ref{suff}.
The existence of the desired structure is, as before, equivalent to
a negative objective of a linear program, and as before, we set up a
slightly simpler linear program first. To wit, the program $L^\prime$
is:

\begin{quote}
\begin{description}

\item[Positivity] All angles are strictly positive.

\item[Euclidean Faces] For every face $t$, $\alpha_t + \beta_t +
  \gamma_t = \pi.$ 

\item[Boundary edge dihedral angles] For every boundary edge $e = AB,$
  incident to the triangle $ABC$, $\gamma + x_e = \Delta(e),$ where
  the slack variables (see the comments following the statement of
  Theorem \ref{duality}) $x_e$ are also non-negative.

\item[Interior edge conditions] For every interior edge $e = AB,$
  incident to  triangles $ABC$ and $ABD$, $\alpha+\delta = \Delta(e).$
\end{description}
\end{quote}

We relax the program $L^\prime$ to a program $L_1^\prime$ by dropping
the requirement that the angles be strictly positive, make the
objective 0, as before, and we see that the dual to the new weakened
linear program $L_1^\prime$ is the following:

\begin{description}

\item[The dual program]
Maximize $F({\bf u}, {\bf v}):\pi \sum_{t \in F(T)} u_t + \sum_{e\in E(T)} \Delta(e)v_e$,
subject to the conditions

\begin{description}
\item[the inequalities of $L_1^*$] $u_t+v_e\leq 0$ whenever $e$ is an
  edge of $t$.

\item[new inequalities] Whenever $e$ is a boundary edge of $T$, $v_e \leq 0$ 
\end{description}

\end{description}

Since the constraints of the above program $L_1^{\prime *}$ are a
superset of the constraints of $L_1^*$, Theorem \ref{maindual} still
tells us that the objective function is maximized if there exists a
$u$, such that $u_t = - v_e = u$, for all $t\in F(T)$, $e\in
E(T)$. Now, this is not enough to guarantee that the 
objective is zero, since the equality case of condition \ref{sumdih}
(when $\mathcal{F} = T$) no longer exists. Since the new inequalities
require $u$ to be non-negative, it follows that for the objective
function to equal $0$, $u$ must be $0$. 

Now, we follow Section \ref{suff} again, to define a program
$L_2^\prime$ in the same way as before (that is, since we want the
angles to be strictly positive, we set $\alpha_i = \beta_i+\epsilon$, etc,
and to also define its dual $L_2^{\prime *}$. By the same reasoning as
before, it follows that the objective of $L_2^{\prime *}$, and hence
of $L_2^\prime$, is negative. In fact, we can do more: we can also
require all of the slack variables $x_e$ to be strictly positive. The
(yet another) new dual program $L_3^{\prime *}$ will have the form:

\medskip\noindent
\textit{Third dual.}
Maximize $F({\bf u}, {\bf v}):\pi \sum_{t \in F(T)} u_t + \sum_{e\in E(T)} \Delta(e)v_e$,
subject to the conditions
\begin{quote}

$u_t+v_e\leq 0$ whenever $e$ is an edge of $t$.

$v_e \leq 0$, when $e$ is a boundary edge of $T$.

$3  \sum_{t \in F(t)} u_t + 
2 \sum_{e\in E(T), e\notin \partial T} v_e + 
2 \sum_{e\in E(T), e \in \partial T} v_e \leq -1$.
\end{quote}

If we omit the last constraint, we remain with the dual program
$L_1^{\prime *}$, and as the discussion above showed, the objective of
that can only be $0$ if $u_t \equiv v_e \equiv 0$, which is at odds
with the last constraint of $L_3^{\prime *}$
}

\thm{theorem}
{
\label{infty}
Let $p_1, \ldots, p_n$ be a set of points in the plane, and $D$ their
Delaunay triangulation. Assume that the shortest edge of $D$ has
length $1$, and the excess of every non-trivial subcomplex of $D$ is
no smaller than $d_0$. Then
$$
\diameter D \leq \left(\frac{4 n}{d_0}\right)^n.
$$
}

The proof of Theorem \ref{infty} will depend on a couple of easy
auxillary results:

\thm{lemma}
{
\label{thmcos}
Let $ABC$ be a triangle with $a/c, a/b < \epsilon$,  $0<\epsilon
<1/10$. Then, $\alpha < 2\epsilon$.
}

\thm{proof}
{
This is an immediate consequence of the theorem of Cosines, or the
theorem of Sines.
}

\thm{lemma}
{
\label{howcom}
Let $\mathcal{F}$ be a subcomplex of $D$. Then the excess of
$\mathcal{F}$ is 
$$
\sum_{\mbox{triangles $ABC$ such that $AB \in \mathcal{F}$}}\gamma(ABC),
$$
where $\gamma(ABC)$ is the angle at $C$.
}

\thm{proof} 
{
This is immediate from the definition of excess.
}

\begin{proof}[of theorem \ref{infty}]
Suppose that the conclusion of the theorem does not hold. Assume,
without loss of generality, that the edge between the vertices $p_0$
and $p_1$ is the shortest one (and thus of length $1$). 

We construct a family of disks $D_1, D_2,\ldots, D_{n}$, all centered
on $p_0$, and such that the radius of $D_i$ is equal to 
$\left(\frac{4 n} d_0\right)^n.$ Let $A_i = D_{i+1} \backslash D_i$,
and let $\mathcal{F}_i$ to be the maximal closed subcomplex of $D$
contained in $D_i$. The hypothesis of the theorem ensures that at
least one of the annuli $A_i$ contains no vertices of $D$, let this
annulus be $A_j$. Then we claim that the excess of $\mathcal{F}_j$ is
smaller than $d_0$. Indeed, consider a triangle $ABC$ of $D$ adjacent to
$\mathcal{F}_j$ along an edge $AB$. The vertex $C$ of $ABC$ lies
outside $D_{j+1}$, and thus the lengths of $AC$ and $BC$ are at least 
$\left(\frac{4n}{d_0}\right)^{j+1} - \left(\frac{4n}{d_0}\right)^j $,
while the length of $AB$ is at most $\left(\frac{4n}{d_0}\right)^j$.
Thus, 
$$\gamma(ABC) < \frac{2}{\frac{4n}{d_0} - 1} < \frac{d_0}{n},$$
by lemma \ref{thmcos}. Thus, by lemma \ref{howcom}, it follows that 
$$\excess \mathcal{F}_j < d_0,$$
contradicting the hypothesis of the theorem.
\end{proof}

\begin{proof}[of Theorem \ref{infmain}]
We only need to show that the positive excess conditions are
sufficient, since they are obviously necessary. In addition, we may
assume that the $1$-skeleton of the Poincar\'e dual $T^*$ is connected
(if not, we prove the theorem for each connected component
separately).

Now, we pick a pair of adjacent base vertices $v_1, v_2 \in V(T),$ and
fix $d(v_1, v_2) = 1$. Now, for $v, w \in V(T)$ we define $d_c(v, w)$ to
be equal to the combinatorial distance between $v$ and $w$ in the
$1$-skeleton of $T$. Now, define $\mathcal{F}_i$ to be the span of all
vertices $u$, such that $d_c(v_1, u) \leq i$. The complex $\mathcal{F}_i$ is 
finite by local finiteness of $T$. In addition, $\bigcup_i
\mathcal{F}_i = T$, so every finite subcomplex of $T$ belongs to some
$\mathcal{F}_i$. For each $\mathcal{F}_i$ we consider a geometric
realization $S(\mathcal{F}_i)$, whose existence is guaranteed by lemma
\ref{linref} (there are many such realizations, we pick any one of them).

Now, enumerate the faces of $T$, in such a way that $t_1$ contains the
edge $v_1 v_2$, and, for any $j$, the faces $t_j$ and $t_{j+1}$ are
adjacent (that is, share an edge) in $T$. For each triangle, we have
the space of shapes (similarity classes), given (for example) by the
complex parameter $z$, obtained by placing the first two vertices of
the triangle at the points $0$ and $1$ in the complex plane, and
reading off the position of the third point (in the upper half-plane,
assuming the triangle is positively oriented). Theorem \ref{infty}
tells us that for any face $t$ of $T$, the set of shape parameters of
realizations of $t$ is contained in a compact set $C_t$ (since the
ratio of lengths of any two sides is bounded by some constant,
depending on the function $\Delta$). We can think of
each $S(\mathcal{F}_i)$ as being an element of $C = \prod_j C_{t_j}$,
which is a compact set by Tykhonov's theorem, and hence we can extract
a convergent subsequence from $S(\mathcal{F}_1), \ldots,
S_(\mathcal{F}_k), \ldots$. Call the limit of that subsequence
$\mathcal{S}$. Since the dihedral angles are obviously continuous
functions of the triangle parameters, the dihedral angles of
$\mathcal{S}$ will be given by $\Delta$, and so $\mathcal{S}$ is the
sought-after realization.
\end{proof}

\medskip\noindent
\textit{Note.} For the comfort of more analytically inclined readers,
it should be pointed out that the last argument is just an
Arzela-Ascoli -- uniform convergence on compact sets argument.

\subsection{Remarks and Corollaries}

Let us assume that the complex $T$ has no boundary, and that the
function $\Delta$ is such that all of the cone angles (computed using
Lemma \ref{cona}) are equal to $2\pi$, so that any realization is
Euclidean (and hence can be developed into the plane). Then, we can
combine theorems \ref{infmain} and \ref{andreevthm} to get 

\thm{corollary}
{
\label{infandreev}
Let $T$ be an infinite locally finite triangulation, and let $\Delta:
E(T) \rightarrow (0, \pi]$ be an assignment of dihedral angles to the
edges of $T$, let Then there exists a flat surface $S$, and a
collection of points $p_1, \ldots, p_n, \ldots$ in $S$ whose Delaunay
triangulation is combinatorially equivalent to $T$ and which has dihedral
angles given by $\Delta$ if and only if:

\begin{itemize}

\item\label{sumvv1} If $v$ is an interior vertex of $T$, then
\begin{equation}
\sum_{\mbox{$e \in E(T)$ incident to $v$}} (\pi - \Delta(e)) = 2\pi.
\end{equation}

\item\label{sumintcyc1} If $\mathcal{E}$ is a non-coterminous minimal
  cutset, then
\begin{equation}
\sum_{e \in \mathcal{E}} (\pi - \Delta(e)) > 2\pi.
\end{equation}

\end{itemize}
}

The above Corollary can be viewed as an extension of Theorem 0.1 of
\cite{steiner} to the case of ideal polyhedra with infinitely vertices
(and hence also of Andre'ev's theorem for ideal polyhedra \cite{Andr70b}),
but not without certain caveats: even when $T$ is topologically a
disk, it is not at all obvious whether the metric on the surface $S$
(or, even, {\em any} surface $S$) is geodesically complete. If it is complete,
it follows that $S$ is the Euclidean plane, and that the developing
map is a global isometry, but otherwise the developing map is,
only an immersion. Checking if the given simply-connected Euclidean
surface $S$ is the Euclidean plane is nothing other than the
\emph{type problem} -- that is, we want to know whether a Riemann surface
is parabolic, hyperbolic, or elliptic. If we know the shapes of all
the triangles, this can be shown to be equivalent to the recurrence of
a random walk on the $1$-skeleton of the complex $T$, where an edge
$AB$, incident to triangles $ABC$ and $ABD$ has weight $1/(\cot C +
\cot D)$ (\cite{irivp}). Using this, it can be shown without too much
difficulty that in the case where all of the dihedral angles are
rational multiples of $\pi$, bounded away from $0$ and $\pi$, then
this is equivalent to the recurrence of the symmetric random walk on
the $1$-skeleton of $T$.  

In the special case of all dihedral angles being equal to $\pi$ or
$\pi/2$, corollary \ref{infandreev} reduces to an existence theorem
for infinite  locally finite circle packings. In this case, the type
problem, and a number of others, has been studied at great length by a
number of authors, starting with Koebe, but more recently by A.~Marden
(in the context of Schottky groups), W.~Thurston, B.~Rodin and
D.~Sullivan, Z.-X.~He, O.~Schramm, and others. 

It should be noted the Theorem \ref{infmain} says nothing about
uniqueness, and the proof certainly does not show any form of
uniqueness. In view of Theorem \ref{uniqueness} one might suspect that
perhaps this could be shown with more work. In fact, it is quite clear
from the above-mentioned work on circle packing that uniqueness fails,
though in a controlled manner described in the conjecture below:

\thm{conjecture}
{
Let $T$ be a infinite locally finite complex, $\Delta$ a system of
dihedral angles satisfying the hypotheses of Theorem \ref{infmain}. If
there is a realization of $T$ supported on a simply-connected domain
$\Omega_0 \subset \mathbf{C}$, then there is one supported on \emph{every}
simply-connected domain $\Omega \subset \mathbf{C}$. Furthermore, such
a realization is determined uniquely by $\Omega$ (up to the group of
M\"obius transformations fixing $\Omega$).
}

\section{Delaunay cells in moduli space}
\label{moduli}
In the Introduction we alluded to the cell decomposition of the moduli
space $\mathcal{M}$ of Euclidean structures, where the cells are given
by Euclidean structures where the Delaunay triangulation has a fixed
combinatorial type $T$. Each cell is a convex polytope
$\mathcal{P}(T)$, and the results above can be used to describe the
combinatorial and geometric structure of $\mathcal{P}(T)$ in some
detail, although some interesting questions (discussed at the end of
this section) remain  open. 

Consider a codimension $1$ face $f$ of $\mathcal{P}(T)$. The face $f$
may be either  a boundary face of the moduli space $\mathcal{M}$
(viewed as a polyhedral complex) or an interior face. In the latter
case, $f$ corresponds to a change of combinatorial type of Delaunay
triangulation from $T$ to $T^\prime$, and it is well understood that
the primitive such change is given by a \emph{diagonal flip}, so $f$
corresponds to a cyclic quadrilateral $ABCD$, where in $T$ the
Delaunay triangulation has triangles (for example) $ABC$ and $CDA$,
while in $T^\prime$, the triangles are $DAB$ and $BCD$. On $f$, the
quadrilateral $ABCD$ can be triangulated either way, but the dihedral
angle along (either) diagonal is equal to $\pi$.

Suppose now that $f$ is a boundary face of $\mathcal{M}$ -- in
particular, this will mean that none of the dihedral angles of the
Euclidean structures in the interior of $f$ are equal to $\pi$. By
Theorems \ref{versimple} and \ref{kappaineq}, it follows that there is
a simple cocycle $c^*$ of $T^*$ where the inequality \ref{kappaineq}[(b)]
becomes an equality. This means, by the geometric estimates in the
beginning of section \ref{ideal},  that the collar $C$ of $c^*$ is
becoming long and thin (that is, the conformal modulus of $C$ diverges
to $\infty$), and so $f$ corresponds to the Euclidean structure
pinching off along the curve $c^*$. The following construction
(discussed in greater detail in an upcoming paper of the author)
helps visualize this pinching off in terms of the more usual
degeneration of hyperbolic surfaces:

Consider a triangulated singular Euclidean surface $(S, E, P)$. 
The data given by $(S, E, P)$ (initially using a triangulation, though
it is easy to show that the result is independent of triangulation)
can be used to construct a cusped hyperbolic surface, as follows:

\begin{enumerate}

\item To each triangle $t$ of the triangulation $T$ of $(S, E, P)$,
  associate an ideal triangle $h(t)$. 

\item For each pair of adjacent triangles $t_1=ABC$ and $t_2=ABD$ of
  $T$, we have the log of the modulus of the cross-ratio of the four
  corresponding points: $$r(t_1, t_2) = \log
  \frac{|AC||BD|}{|BC||AD|}.$$

\item For each pair of adjacent triangles $t_1$ and $t_2$ as above,
  glue the hyperbolic triangles $h(t_1)$ and $h(t_2)$ along the 
  edge corresponding to $AB$ with shear   (see, eg, \cite{rivid})
  equal to $r(t_1, t_2).$ 

\end{enumerate}

It is not hard to see that if we start with a Euclidean structure $(S,
E, P)$, we will wind up with a \emph{complete} cusped hyperbolic
structure (actually, it is sufficient, but not necessary, to start
with a Euclidean structure -- some similarity structures will give a
complete structure also). The construction thus defines a map
(certainly not injective, but which can be shown to be surjective
using the construction of \cite{penner,ep,np})
between the moduli space of Euclidean structures with cone points and
that of complete finite-area hyperbolic structures. It can be shown
(this was an important part of \cite{rh,steiner} for the case of genus
$0$) that the degeneration, as above, of the Euclidean structure on
$(S, E, P)$, corresponds precisely to the pinching off along a simple
closed curve of the hyperbolic structure on $(S, H(E), P)$.

\medskip\noindent
\textit{Remark.} Another surjection between the space of singular
Euclidean structures and the space of cusped hyperbolic structures is
well-known, and is discussed in the well-known paper of Troyanov
\cite{troy}. The two surjections are \emph{not} the same -- this was
remarked by C.~T.~McMullen.

\section{Computational complexity}
\label{complexsec}

Consider the following two decision problems:

\begin{description}
\item[Problem 1]
Let $\mathcal{T}$ be a simplicial complex, homeomorphic to a surface $S$
(possibly with boundary), and let $\alpha: V(\mathcal{T}) \rightarrow
{\bf R}^+$ be an assignment of cone angles to the vertices of
$\mathcal{T}$. Does there exist a Euclidean structure $E$ on $S$ with
the prescribed cone angles, such that the Delaunay triangulation of $(S,
E)$ is combinatorially equivalent to $\mathcal{T}$?

\item[Problem 2]
Let $\mathcal{T}$ be a simplicial complex, homeomorphic to a surface
$S$ (possibly with boundary), and let $\Delta: E(\mathcal{T})
\rightarrow (0, \pi]$ be an assignment of dihedral angles to the edges
of $\mathcal{T}$. Does there exist a singular Euclidean structure $E$
on $S$, such that the Delaunay triangulation of $(S, E)$ is
combinatorially equivalent to $\mathcal{T}$, and whose dihedral angles
are given by $\Delta$.
\end{description}

By Theorem \ref{exist}, there are efficient algorithms for both
problems, since they reduce to the linear program $L$ of section
\ref{oldstuff}. If the angles (cone angles in problem 1, dihedral
angles in problem 2) are rational multiples of $\pi$, such that the
numerator and denominator are both bounded by $C$, then (by now)
standard interior point methods allow us to solve the linear program
$L$ using $O(n^4 (1+ \log C))$ arithmetic operation, each involving
arithmetic using precision $O(n (1+\log C))$. In the case where all
of the prescribed cone angles are equal to $2\pi$, the $\log C$ can be
disposed with, and we wind up with an algorithm of bit-complexity
$O(n^5 \log^2 n).$ 

\medskip\noindent
\textit{Remark.} In practice, the simplex algorithm appears much more
efficient, and this has been used by M.~B.~Dillencourt to 
analyze all planar triangulations of up to $14$ vertices, and to
determine which of them are combinatorially equivalent to planar
Delaunay tessellations (\cite{Dill1,Dillp}).

For problem 2, the network flow formulation of section \ref{nflow}
turns out to give a markedly superior complexity. Indeed, it has been
shown in \cite{aot} that for a network with $n$ nodes,
$m$ arcs, and the (integer) capacity of each arc bounded by $U$, we
can determine the maximal flow in time bounded by $O(n m
\log((n/m)(\log U)^{1/2}+2))$. 

For the network $N_1$ of section
\ref{nflow}, assuming the genus of the surface is fixed, the number of
arcs and nodes in the network are both bounded by constant multiples
of the number $F$ of faces in the complex $T$. If all of the dihedral
angles are rational multiples of $\pi$, with numerators and
denominators all bounded in absolute value by $C$, the quantity $U$ is
bounded by $C^{O(F)}$ (since we need to compute the least common
multiple of the denominators), giving a running time bound of
$O(F^{5/2}).$ For the program $N_2$, the bound is the same, since the
only difficulty consists of picking the right value of $\epsilon$, and
this can be made to be $1/\mbox{least common denominator of the
dihedral angles}$.

\section{Some applications to combinatorial geometry}
\label{someapps}

A well-known theorem of Steinitz says that every three-connected 
planar graph can be realized as the $1$-skeleton of a convex 
polyhedron in $R^{3}$, while a famous theorem of Aleksandrov states 
that every Euclidean metric on $S^{2}$ with positively curved 
cone-points can be realized as the induced metric on the surface of a 
convex polyhedron in $R^{3}$. Below, we show a negative result, which 
should be compared with the final example of \cite{rit}:

\thm{theorem}
{
\label{badtri}
There exist infinitely many triangulations of $S^{2}$ which can not be
realized as a  Delaunay triangulation with respect to the cone-points
of any Euclidean  metric on $S^{2}$ with positively curved cone-points.
}

First a definition:

\thm{definition}
{
Let $T$ be a triangulation. The \emph{stellation} $s(T)$ of $T$ is the 
complex obtained by replacing each face $ABC$ of $T$ by three faces 
$AOB$, $AOC$, and $BOC$.
}

Theorem \ref{badtri} follows immediately from the following claim.

\thm{claim}
{
Let $T$ be any triangulation of $S^{2}$ with at least eight faces. 
Then the stellation $s(T)$ of $T$ is not combinatorially equivalent 
to the Delaunay triangulation of any Euclidean metric on $S^{2}$ with 
positively curved cone-points, where the cone-points correspond to 
vertices of $s(T)$.
}

\begin{proof} Let an \emph{old} vertex of $s(T)$ be one that 
was already a vertex of $T$, while a \emph{new} vertex be one that was 
added at stellation. The set $\mathcal{N}$ of new vertices of $s(T)$ corresponds 
to the  set of  faces of $T$. For any vertex $v$, recall that $C(v)$ 
denotes the cone angle at $v$. The Gauss-Bonnet theorem tells us that
\begin{equation}
\sum_{v\in V(s(T))}(2\pi - C(v)) = 4\pi.
\end{equation}
Or, recombining the terms:
\begin{equation}
\label{gaussb}
\sum_{v\in V(s(T))}C(v) = 2\pi (|V(s(T))| - 2).
\end{equation}
Note that every edge of $s(T)$ is incident to an old vertex, and to 
at most one new vertex. Combining this observation with Lemma \ref{cona},
we see that for any \emph{Delaunay} triangulation combinatorially
equivalent to $s(T)$, it must be true  that
\begin{equation}
\label{badineq}
\sum_{v \in \{\mbox{old vertices of $s(T)$}\}} C(v) \geq 
\sum_{v \in \{\mbox{new vertices of $s(T)$}\}} C(v).
\end{equation}
Combining Eq. \ref{badineq} with Eq. \ref{gaussb}, it follows that
\begin{equation}
\label{badderineq}
\sum_{v \in \{\mbox{old vertices of $s(T)$}\}} C(v) \geq \pi (|V(s(T)) - 
2).
\end{equation}
Note now that $|V(s(T))| = |V(T)| + |F(T)|$. By a standard 
computation using Euler's formula for triangulations of the sphere, 
we know that $|V(T)| = \frac12 |F(T)| + 2,$ thus $|F(T)| = 2 |V(T)| - 4.$
Thus, 
\begin{equation}
\label{neww}
|V(s(T))-2| = 3 |V(T)| - 6.
\end{equation}
Equations \ref{neww} and \ref{badderineq} together imply that the
\emph{average} cone angle at an old vertex of $s(T)$ is at least 
$\pi(3 |V(T)| - 6)/|V(T)| = 3 \pi - 6/|V(T)|.$
If $|V(T)| > 6$ it follows that the average cone angle at an old
vertex is greater  than $2\pi$, which contradicts the assumption that
the cone angles were positively curved.
\end{proof}

\section{Linear hyperbolic structures on $3$-manifolds}
\label{lhyp}

As explained in \cite{rvol}, the study of Euclidean triangulations on
surfaces is essentially equivalent to the study of hyperbolic ideally
triangulated complexes, which are combinatorially just cones over the
triangulated surface. Hence, the contents of this section are closely
related to the subject-matter of much of the rest of the paper, in
more than just the linear programming approach.

Consider a $3$-manifold $M^3$ with boundary a collection of tori, and
consider a topological ideal triangulation $T$ of $M^3$. We would like
to know when there is a complete hyperbolic structure on $M^3$, such
that $T$ is a geometric ideal triangulation. In general, this is a
very difficult question, at least as hard as Thurston's
hyperbolization conjecture (since even if $M^3$ admits a hyperbolic
structure of finite volume, there might not be an ideal triangulation
combinatorially equivalent to $T$). However, below we will consider a
``linear'' version of the question above. 

Recall that an ideal simplex $S$ in $\mathbf{H}^3$ has the properties that

\begin{description}

\item[Euclidean links]
The link of each vertex is a Euclidean triangle.

\item[Equal opposite dihedral angle]
If $S = ABCD$, then the dihedral angles at the edges $AB$ and $CD$ are
equal (this is actually a consequence of the condition on the links).
\end{description}

An ideal simplex is thus determined by the angles of the link of any
one of its vertices (all links are easily seen to be the same).

Now, if $T$ comes from a genuine hyperbolic structure, it must be true
that the sum of the dihedral angles incident to the edges of $T$ must
equal $2\pi$, and so for $T$ to correspond to such a structure, the
following linear program must have a strictly positive solution:

\begin{description}
\item[The variables] These are the dihedral angles of the simplices. For
  each simplex $S$ we use three angles $\alpha, \beta, \gamma$
  corresponding to the angles of the link of one vertices of $S$.

\item[Simplex conditions]
For each simplex $S$, the sum of the dihedral angles is equal to $\pi$:
\begin{equation}
f_S: \alpha + \beta + \gamma = \pi.
\end{equation}

\item[Edge conditions]
For each edge $e$ of $T$ the sum of the dihedral angles of all
simplices incident to $e$ equals $2\pi:$
\begin{equation}
f_e: \sum \lambda_i = 2\pi.
\end{equation}
\end{description}

\thm{definition}
{
We say that if the above linear program has a non-negative solution,
then the pair $(M^3, T)$ admits a \emph{weak linear hyperbolic
  structure}. If the linear program has a strictly positive solution,
then we say that the pair $(M^3, T)$ admits a \emph{linear hyperbolic
  structure.}
}

It is clear that the existence of a linear hyperbolic structure is, in
general, no guarantee of the existence of a genuine hyperbolic
structure, since the linear conditions do not preclude ``Dehn
surgery'' singularities, as well as translational singularities along
edges of $T$ (see \cite{thgt3m,nz} for discussion). Conversely, as
remarked above, the existence of a 
complete hyperbolic structure on $M^3$ is no guarantee that there is a
positively oriented ideal triangulation combinatorially equivalent to
$T$ (or, indeed, any positively oriented triangulation). However,
linear hyperbolic structures have the advantage of being considerably
more tractble. 

\thm{theorem}
{
\label{hext}
In order for there to exist a weak linear hyperbolic structure for
$(M^3, T)$, every normal surface with respect to $T$ must have
non-negative Euler characteristic. In order for there to exist a
linear hyperbolic structure for $(M^3, T)$, every
non-boundary-parallel normal surface with respect to $T$ must have
strictly negative Euler characteristic.
}

\thm{proof}
{
We will use the method of Section \ref{suff}. First, let us write the
dual program of the linear program (referred to as $L_h$ in the
sequel) for weak hyperbolic structure: Our
variables are $\{v_S\}$, where $S$ ranges over all the simplices of
$T$, and $\{v_e\}$ where $S$ ranges over all the edges of $T$. 

The dual program $L_h^*$ is: 

\begin{description}

\item{Maximize} $\sum_{S\in S(T)} v_S + 2 \sum_{e \in E(T)} v_e$.

\item{Subject to} $v_S + v_{e_1} + v_{e_2}\leq 0$ for all faces $S$
  and pairs of opposite edges $e_1, e_2$.
\end{description}

In order for the primal program to have a non-empty feasible region,
the objective of the dual must be non-positive.

Consider now a normal surface $\mathcal{S}$ (see, eg,
\cite{hack,jaco,hemion} for rudiments of normal surface theory). The
surface $\mathcal{S}$ intersects each simplex $S$ in a collection of
disks, which are combinatorially either triangles (cutting off one
vertex of $S$ from the other three), or quadrilaterals (separating one
pair of vertices from another). For each simplex $S$, define
$t_S(\mathcal{S})$ to be the number of triangular components of
$\mathcal{S} \cap S$, and $q_S(\mathcal{S})$ to be the number of
quadrilateral components. For each edge $e$ 
of $T$, define $i_e(\mathcal{S})$ to be the number of intersections of
$\mathcal{S}$ with $e$. The intersections of $\mathcal{S}$ with the
simplices of $T$ induces a triangulation $\mathcal{T}$ of $S$, where each triangular
disk contributes one triangle, and each quadrangle contributes
two. Define $u_S(\mathcal{S})$ to be the number of triangles of
$\mathcal{T}$ sitting inside a simplex $S$. Evidently, $u_S = t_S + 2
q_S$. 

Note that, by Euler's formula, $\chi(\mathcal{S}) = \sum_{e \in E(T)}
i_e -\frac12 \sum_{S \in S(T)} u_S$. This is seen to be very similar in form
to the objective function of $L_h^*$, so let us set
$v_S = - u_S$, and $v_e = i_e$. 

\thm{lemma}
{
\label{normvar}
The assignment of the variables as above satisfies the inequality
constraints of $L_h^*$.
}

\begin{proof}[of lemma]
We need to check that for $S$ a simplex and $e_1, e_2$ a pair of
disjoint edges of $S$. We need to check that 
\begin{equation}
\label{normineq}
i_{e_1} + i_{e_2} \leq u_S.
\end{equation}
By linearity, we need just check the inequality \ref{normineq} for
connected components of $\mathcal{S} \cap S$. If that component is a
triangle $t$, then $t$ contributes $1$ to $i_{e_1} + i_{e_2}$ (since a
``normal triangle'' intersects exactly one of each pair of opposite
edge). Also, $t$ contributes $1$ to $u_S$, so for a triangular face,
the right and left hand sides of (\ref{normineq}) are equal.

Suppose now we have a quadrilateral component $q$. It 
contributes $2$ to the right hand side of (\ref{normineq}). As for the
left hand sides, $q$ hits two pairs of opposite sides of $S$, so if
$e_1$ and $e_2$ is one of those pairs, then we have a contribution of
$2$ to the left hand side, and otherwise we have a contribution of
$0$.
\end{proof}

\thm{remark}
{
\label{rmkkk}
Notice that if $\mathcal{S}$ is such that all of the
components of $\mathcal{S} \cap S$, for all $S\in S(T)$ are triangles,
then all of the constraints of $L_h^*$ are equalities with the
assignment of variables as above. Any such $\mathcal{S}$ is easily
seen to be a union of boundary tori.
}

Lemma \ref{normvar} concludes the proof of the ``weak'' part of
Theorem \ref{hext}, since if any $\mathcal{S}$ had positive Euler
characteristic, the program $L_h^*$ would have a positive objective,
and thus the program $L_h$ would have no solution.

For the proof of the ``strong part'' we use the same trick as in
Section \ref{suff}. Define new variables $\alpha^\prime = \alpha +
\epsilon$, \textit{etc.}
Our primal linear hyperbolicity program $L_s$ is now:

\begin{description}

\item{Minimize:} $-\epsilon$

\item{subject to face constraints} $\alpha^\prime + \beta^\prime +
  \gamma^\prime + 3\epsilon = \pi.$

\item{and to edge constraints} $\sum_{\alpha^prime} + v(e) \epsilon =
  2\pi$, where $v(e)$ is the valence of $e$.
\end{description}

The dual program $L_s^*$ is then

\begin{description}
\item{Maximize} $\sum_{S\in S(T)} v_S + 2 \sum_{e \in E(T)} v_e$.

\item{Subject to the old constraints} $v_S + v_{e_1} + v_{e_2}\leq 0$
  for all faces $S$ and pairs of opposite edges $e_1, e_2$.

\item{and the new constraint} $3\sum_{S \in S(T)} v_S + \sum_{e\in
    E(T)} v(e) v_E \leq -1.$
\end{description}

In order for $(M^3, T)$ to be linearly hyperbolic, the objective must
be strictly negative.

Observe that the sum of the left hand sides of the old constraints is
\emph{equal} to the left hand side of the new constraint. Indeed, each
$v_S$ occurs three times (once for each pair of opposite edges), and
each $v_e$ occurs the number of times equal to the valence of
$e$. Hence, the new constraints simply says that in at least one of
the old constraints the inequality must be \emph{strict}. Keeping in
mind Remark \ref{rmkkk}, this implies that every non-boundary-parallel
normal surface must have strictly negative Euler characteristic, thus
proving the second part of Theorem \ref{hext}
}

Some remarks may be in order: 
It is easy to see (and not surprising) that an identical theorem can
be proved if the cone angles around the edges of $T$ are required to
not be smaller than $2\pi$, while if the angles are smaller than
$2\pi$, one can show an analogous ``orbifold'' version of the theorem.
A more interesting question is whether theconverse of Theorem
\ref{hext} holds. This is equivalent to asking whether every
assignment of variables satisfying the constraints of programs $L_h^*$
and $L_s^*$ comes from a normal surface. 

\begin{acks}
The author would like to thank the \'Ecole Normale Superieure de Lyon
and the Institut Henri Poincar\'e for their hospitality. The special
case of Theorem \ref{mainthm} for Euclidean structures (that is,
trivial holonomy) has been independently shown for closed surfaces by
W.~Veech (\cite{veech}) and B.~H.~Bowditch (\cite{bhb}) using
completely different methods. Theorem \ref{hext} is due to A.~Casson. 
\end{acks}

\bibliographystyle{acmtrans}

\begin{thebibliography}{99}

\bibitem{aot}
K.~Ahuja,  J.~B.~Orlin, and R.~E.~Tarjan.
\textit{Improved time bounds for the maximum flow problem}, SIAM
Journal on Computing,
\textbf{18}(1989), pp. 939-954.

\bibitem{Andr70b}
E.~M. Andreev.
\textit{On convex polyhedra of finite volume in Lobachevskii space.\/}
Math. USSR, Sbornik, \textbf{12}:255--259, 1970.

\bibitem{Dillp}
M.~B.~Dillencourt, private communication.

\bibitem{Dill1}
M.~Dillencourt. 
\textit{Polyhedra of small order and their Hamiltonian properties\/}
J. Combin. Theory Ser. B \textbf{66}(1), 1996, pp. 87--122.

\bibitem{be}
B.~H.~Bowditch and D.~B.~A.~Epstein. 
\textit{Natural triangulations associated to a surface.\/}
Topology \textbf{27}(1), pp. 91-117.

\bibitem{bhb}
B.~H.~Bowditch.
\textit{Singular Euclidean structures on surfaces.\/} Journal of the
London Mathematical Society \textbf{44}(3), pp 553-565.

\bibitem{ep}
D.~B.~A.~Epstein and R.~C.~Penner.
\textit{Euclidean decompositions of non-compact hyperbolic manifolds},
Journal of Differential Geometry \textbf{27} (1988), pp 67-80.

\bibitem{hack}
W.~Haken.
\textit{Theorie der Normalfl\"achen}, (German) Acta Math. \textbf{105}, 1961,
pp. 245--375.

\bibitem{hemion}
G.~Hemion.
The classification of knots and $3$-dimensional spaces. Oxford
Science Publications. The Clarendon Press, Oxford University Press,
New York, 1992.

\bibitem{harer}
John Harer.
\textit{The virtual cohomological dimension of the mapping class group
  of an orientable surface.\/} 
Invent. Math. \textbf{84} (1986), no. 1, pp.~157--176. 

\bibitem{hrs:ann}
C.~D.~Hodgson, Igor Rivin and W.~D.~Smith.
\textit{A characterization of convex hyperbolic polyhedra and
of convex polyhedra inscribed in the sphere\/},
Bulletin of the American Mathematical Society \textbf{27}(3), October 1992. 

\bibitem{jaco}
W.~Jaco. 
Lectures on three-manifold topology. CBMS Regional Conference Series in
Mathematics, \textbf{43}. American Mathematical Society, Providence,
R.I., 1980. 

\bibitem{np}
M.~N\"a\"at\"anen and R.~C.~Penner.
\textit{The convex hull construction for compact surfaces and the
  Dirichlet polygon}, Bulletin of the London Mathematical Society
\textbf{6}, 1991. pp 568-574.

\bibitem{nz}
W.~D.~Neumann and D.~Zagier.
\textit{Volumes of hyperbolic $3$-manifolds}, Topology \textbf{24}(3),
pp. 307-332.

\bibitem{penner}
R.~C.~Penner.
\textit{The decorated Teichm\"uller space of punctured surfaces,\/},
Comm. Math. Phys. {\bf 113}(2), 1987, pp.~299-339,

\bibitem{rh}
Igor Rivin and C.~D.~Hodgson.
\textit{A characterization of compact convex polyhedra in hyperbolic
  $3$-space\/}, Inventiones Mathematicae \textbf{111}(1), 1993.

\bibitem{rit}
Igor Rivin.
\textit{On geometry of convex ideal polyhedra in hyperbolic $3$-space\/},
Topology \textbf{32}(1), 1993.

\bibitem{rivid}
Igor Rivin.
\textit{Intrinsic geometry of convex ideal polyhedra in hyperbolic
  $3$-space\/}, 
Analysis, algebra, and computers in mathematical research (Proceedings
of the 1992 Nordic Congress of Mathematicians, Luleaa, 1992), Lecture
Notes in Pure and Applied Mathematics 156, Marcel Dekker, 1994.

\bibitem{rvol}
Igor Rivin. 
\textit{Euclidean structures on simplicial surfaces and hyperbolic
volume\/}, Annals of mathematics \textbf{139}(3), 1994.

\bibitem{steiner}
Igor Rivin.
\textit{A characterization of ideal polyhedra in hyperbolic 3-space\/},
Annals of mathematics \textbf{143}(1), 1996.

\bibitem{irivp}
Igor Rivin. Personal communication.

\bibitem{st:geom}
Jakob Steiner,
Systematische Entwicklung der Abh\"angigkeit
geometrischer   Gestalten von einander,
Reimer, Berlin, 1832; Appeared in J. Steiner's Collected Works, 1881.

\bibitem{Stei1}
E.~Steinitz.
\textit{Uber isoperimetrische Probleme bei konvexen Polyedern\/},
J.~Reine Angew.~Math. \textbf{159}, pp.~133-143, 1928.

\bibitem{troy}
M.~Troyanov.
\textit{Les surfaces euclidiennes \`a singularites coniques},
Enseignement Mathematique, \textbf{32}(1-2), 1986, pp. 79-94.

\bibitem{thgt3m}
W.~P.~Thurston.
\textit{Geometry and topology of $3$-manifolds\/}, Unpublished
Princeton lecture notes for 1978, distributed by MSRI, Berkeley.

\bibitem{vlw}
J.~van~Lint, R.~M.~Wilson.
A Course in Combinatorics, Cambridge University Press, Cambridge,
1992.

\bibitem{veech}
W.~Veech.
\textit{Delaunay Partitions\/},
Topology {\bf 36}, no.~1, 1997, pp.~1-28.

\end{thebibliography}

\end{document}